\documentclass[12pt,a4paper]{article}

\usepackage{amsmath,amssymb,amsthm}
\usepackage[shortlabels]{enumitem}
\setlist[enumerate]{label=(\roman*)}
\usepackage{graphicx}
\usepackage{tikz}
\usepackage[colorlinks=false]{hyperref}
\hypersetup{
    colorlinks,
    linkcolor={red!50!black},
    citecolor={blue!50!black},
    urlcolor={blue!80!black}
}
\usepackage[capitalise]{cleveref}
\usepackage{xcolor}
\usepackage{calrsfs}
\usepackage{natbib}
\setcitestyle{semicolon}

\newtheorem{theorem}{Theorem}

\newtheorem{cor}[theorem]{Corollary}
\theoremstyle{definition}

\newtheorem{conj}{Conjecture}

\newcommand \dm[1]  { \,\mathrm d{#1} }

\renewcommand{\epsilon}{\varepsilon}
\renewcommand{\phi}{\varphi}

\renewcommand{\leq}{\leqslant}
\renewcommand{\geq}{\geqslant}

\newcommand{\norm}[1]{\left\lVert#1\right\rVert}

\newcommand{\norml}[3]{\norm{#1}_{L^{#2}({#3})}}

\newcommand{\pair}[1]{\left\langle{#1}\right\rangle}

\newcommand{\abs}[1]{\left\vert#1\right\vert}

\renewcommand{\hat}[1]{\widehat{#1}}

\newcommand{\e}{\mathrm{e}}

\newcommand{\tpi}{2 \pi \mathrm{i}}

\newcommand{\R}{\mathbb{R}}
\newcommand{\Rn}{\R^n}
\newcommand{\C}{\mathbb{C}}
\newcommand{\N}{\mathbb{N}}

\newcommand{\Z}{\mathbb{Z}}
\newcommand{\Zn}{\Z^n}
\newcommand{\T}{\mathbb{T}}
\newcommand{\Tn}{\T^n}

\newcommand{\D}{\mathbb{D}}

\renewcommand{\th}{\textsuperscript{th} }
\newcommand{\qtext}[1]{\quad\text{#1}\quad}
\newcommand{\qand}{\qtext{and}}

\title{Remarks on the $L^p$ convergence of Bessel--Fourier series on the disc}
\author{Ryan L. Acosta Babb%
\thanks{University of Warwick, UK
(\href{mailto:r.acosta-babb@warwick.ac.uk}{r.acosta-babb@warwick.ac.uk}).}}
\date{}

\begin{document}
\maketitle

\begin{abstract}
    The $L^p$ convergence of eigenfunction expansions for the Laplacian on planar domains
    is largely unknown for $p\neq 2$.
    After discussing the classical Fourier series on the 2-torus, we move onto the disc,
    whose eigenfunctions are explicitly computable as products of trigonometric and Bessel functions.
    We summarise a result of \cite{CordobaBalodis} regarding the $L^p$ convergence
    of the Bessel--Fourier series in the mixed norm space $L^p_{\mathrm{rad}}(L^2_{\mathrm{ang}})$
    on the disk for the range $\tfrac{4}{3}<p<4$. We then describe how to modify their result
    to obtain $L^p(\D, r\dm{r}\dm{t})$ norm convergence in the subspace $L^p_{\mathrm{rad}}(L^q_{\mathrm{ang}})$
    ($\tfrac{1}{p}+\tfrac{1}{q}=1$) for the restricted range ${2\leq p < 4}$.
\end{abstract}
%\tableofcontents

\section{Introduction}
For a function $f\in L^2(\Tn)$, we can truncate its Fourier series by ``spherical modes''
\begin{equation}\label{eqn:sphmodes}
    S_{N}f := \sum_{\abs{k}\leq N}\hat{f}(k)\e^{\tpi k\cdot x},
\end{equation}
or by ``cubic modes''
\begin{equation}\label{eqn:cubmodes}
    S_{[N]}f := \sum_{\abs{k_j}\leq N}\hat{f}(k)\e^{\tpi k\cdot x},
\end{equation} where \[
    k = (k_1,\ldots,k_n)\in\Zn \qand \abs{k}^2 = \sum_{j=1}^n\abs{k_j}^2.
\]
It is well known that $S_Nf$ from \cref{eqn:sphmodes} fails in general to converge to $f$ when $p\neq 2$.
This follows, by standard transference arguments \citep[see][]{GrafakosCFA},
from Fefferman's (\citeyear{Fefferman1971}) result that the
indicator function of the ball is not an $L^p$-bounded Fourier multiplier for any $p\neq 2$.
\citep[See][for a detailed discussion and related results.]{Grafakos2014MFA}
On the other hand, the square truncations from \cref{eqn:cubmodes} are perfectly well-behaved
for all $1<p<\infty$ \citep[see again][]{GrafakosCFA}.

Recently, \cite{RobinsonFefferman2021} have asked whether,
given a differential operator with an orthonormal family $w_k$ of eigenfunctions,
there is a choice $\{w\in E_N : N\in\N\}$ of eigenfunctions such that the ``truncations'' \[
    S_{E_N}f = \sum_{w\in E_N}\pair{f,w}w
\] are ``well-behaved'' in $L^p$ for all $1<p<\infty$.

That this is not possible in general can be shown by considering the disc $\D\subset\R^2$.
The eigenfunctions for the Laplacian on $\D$ are of the form \[
    \e^{\tpi \theta m}J_{m}(j_m^nr)
    \qtext{for} (r,\theta)\in[0,1]^2, (m,n)\in\Z\times\N
\] corresponding to the respective eigenvalues $4\pi^2m^2+(j_n^m)^2$.
Here $J_m:=J_{\abs{m}}$ denotes a Bessel function of the first kind and $j_n^m:=j_n^{\abs{m}}$ its non-negative zeros 
\citep[see][]{Watson1995}.

Consider the function $f(r)=r^{-3/2}$, which lies in the space ${L^p([0,1],r\dm{r})}$ for $1\leq p < 4/3$.
\cite{Wing1950} proved that, for any choice of $J_m$, the 1-dimensional Bessel series of $f$ fails to converge in $L^p([0,1],r\dm{r})$.
By letting $g(r,t) := f(r)$ in $L^p(\D)$, it follows that the 2-dimensional Bessel--Fourier series of $g$ is \[
    \sum_{m\in\Z}\sum_{n\in\N} a_{m,n}J_m(j_n^mr)\e^{\tpi mt} = \sum_{n\in\N}a_nJ_m(j_n^mr)
\] and so does not converge to $g$ for any $m\geq 0$, \emph{regardless of the truncation method}.

Thus, restrictions on the range of $p$ are to be expected.
A natural range is $4/3 < p <4$, since this is precisely the range that works for the 1-dimensional Bessel series
\citep{Wing1950}.
\citep[It is instructive to compare this to the ranges of $L^p$ convergence for the Bochner-Riesz means on $\R^2$; see]%
[for details. We will return to this question later.]{GrafakosCFA}

The best result known at this time is due to \cite{CordobaBalodis}, who reduced the problem of convergence on the disc
to extant results on the convergence of Fourier and 1-dimensional Bessel series, albeit with a modified norm.
We will exploit their argument to obtain $L^p$ convergence in a certain subspace of $L^p$.
\section{Mixed-norm convergence}
Define the space $L^p_{\mathrm{rad}}(L^2_{\mathrm{ang}})$ by the inequality
\begin{equation}\label{eqn:p2norm}
    \norm{f}_{p,2} := \left[\int_0^1\left(\sum_{m}\abs{f_m(r)}^2\right)^{p/2}r\dm{r}\right]^{1/p}
    \equiv \norml{\norm{(f_m(r))}_{\ell^2}}{p}{r\dm{r}} < \infty,
\end{equation}
where $f_m(r)$ are the Fourier coefficients of the angular function $t\mapsto f(r,t)$: \[
    f_m(r) := \int_{0}^1 f(r,\theta)\e^{-\tpi m\theta}\dm{\theta}.
\]
Denote by $S_{N,M}f$ the partial sums of the Bessel--Fourier series of $f:\D\to\C$:
\begin{equation}
    S_{N,M}f(r,t) := \sum_{m=-M}^{M}\sum_{n=1}^{N}a_{m,n}J_{m}(j_n^mr)\e^{\tpi mt}.
\end{equation}
(We drop the superscript ``$(d)$" present in \citealp{CordobaBalodis}, since $d=2$ will remain fixed in our discussion.)
\begin{theorem}[\citealp{CordobaBalodis}]\label{thm:CBUBT}
    The operators $S_{N,M}$ are uniformly bounded on $L^p_{\mathrm{rad}}(L^2_{\mathrm{ang}})$
    if, and only if, $\frac{4}{3}<p<4$ when $N\geq AM+1$ for an absolute constant $A>0$.
\end{theorem}
The norm convergence of the series to $f$ follows by the usual uniform boundedness argument.
(See \citealp{CordobaBalodis};
cf. the analogous Fourier series argument in \citealp{GrafakosCFA}.)

To attack the proof of \cref{thm:CBUBT}, they exploited the presence of the Fourier coefficients
in the norm (\ref{eqn:p2norm}).
Indeed, $S_{N,M}f$ is a trigonometric polynomial whose $m$\th Fourier mode ($\abs{m}\leq M$) is
\begin{equation}\label{eqn:1dfbop}
    S_{N,m}f_m(r) \equiv \sum_{n=1}^N a_{m,n}J_m(j_n^mr),
\end{equation}
which is precisely the 1-dimensional Bessel series summation operator for the radial function $r\mapsto f_m(r)$
in terms of the $m$\th order Bessel function $J_m$.
Thus,
\begin{align*}
    \norm{S_{N,M}f}_{p,2} &=
    \left[  \int_0^1
                \left( \sum_{m=-M}^M \abs{S_{N,m}f_m(r)}^2 \right)^{p/2}
            r\dm{r}   \right]^{1/p}\\
    &= \norml{(S_{N,m} f_m)_m}{p}{r\dm{r} ; \ell^2},
\end{align*}
so the boundedness of $S_{N,M}$ in $L^p_{\mathrm{rad}}(L^2_{\mathrm{ang}})$ is reduced
to a uniform bound for vector-valued inequalities.
Such bounds must be independent of the length, $2M+1$, of the vector \[
    \left(S_{N,m}f_{-m},\ldots,S_{N,m}f_m\right).
\] Note that $S_{N,-m}=S_{N,m}$ by our convention that $J_{m}=J_{\abs{m}}$ for $m\in\Z$.
The functions $f_m$ and $f_{-m}$, however, are distinct in general,
as they correspond to distinct Fourier coefficients.

Let us now turn to $L^p$ convergence on the disc, where the relevant norm is \[
    \norml{f}{p}{\D} = \norml{\norml{f(r,t)}{p}{\dm{t}}}{p}{r\dm{r}}.
\]
For $p\neq2$ we cannot replace the inner ``angular'' $L^p$ norm by a sum of Fourier coefficients
\citep[see Chap. IV of][]{Katznelson2004}.
However, for $p\geq2$ we have the \emph{Reverse Hausdorff-Young Inequality}:
\begin{equation*}%\label{prop:revHausYoung}
        \norml{g}{p}{\T} \leq \norm{(\hat{g}(k))_k}_{\ell^q(\Z)}
        \qtext{whenever}  p \geq 2 \text{ and } \frac{1}{p}+\frac{1}{q}=1.
\end{equation*}
We therefore have
\begin{align*}%\label{eqn:LpLpqineq}
    \norml{f}{p}{\D} &= \left[\int_0^1 \norml{f(r,t)}{p}{\T,\dm{t}}^pr\dm{r}\right]^{1/p}\nonumber\\
    &\leq \left[\int_0^1\left(\sum_k\abs{f_k(r)}^q\right)^{p/q}r\dm{r}\right]^{1/p}\nonumber\\
    &=: \norm{f}_{p,q}.
\end{align*}
Using this norm, we define the space \[
    L^p_{\mathrm{rad}}(\ell_{\mathrm{ang}}^q) := \left\{
        f\in L^p(\D) : \norm{f}_{p,q} <\infty
    \right\}  .
\]

Careful inspection of the proofs in \cite{CordobaBalodis} shows that the space
$\ell^2$ can be replaced by $\ell^q$ throughout when \[
    \frac{1}{p} + \frac{1}{q} = 1 \qand 2 \leq p < 4.
\] The kernel of the 1-dimensional summation operators $S_{N,m}$ in \cref{eqn:1dfbop}
are controlled by weighted, vector-valued norms on the operators
\begin{equation}\label{eqn:niceops}
    \int_{0}^1\frac{f(t)}{2-x-t}\dm{t}\qand
    \int_0^1\frac{f(t)}{x+t}\dm{t},
\end{equation} the Hilbert Transform and the Hardy-Littlewood Maximal Functional.
The weight $r^{1-p/2}$ satisfies the 1-dimensional Muckenhoupt $A_p$ condition if, and only if,
\[
    -1 < 1-\frac{p}{2} < p-1 \qtext{that is} \frac{4}{3} < p < 4
\]
\citep[see][Example 7.1.7]{GrafakosCFA} and, when this is the case,
we have the inequalities \[
    \norml{(\mathcal{M}f_k)_k}{p}{r^{1-p/2};\ell^q}
    + \norml{(Hf_k)_k}{p}{r^{1-p/2};\ell^q}
    \lesssim_{p,q} \norml{(f_k)_k}{p}{r^{1-p/2};\ell^q}.
\]\cite[See][(B) on p.\,25.]{Cordoba1989}
Furthermore, the kernels of the operators in \cref{eqn:niceops} are nice enough
that both are bounded on the space $L^p([0,1],r^{1-p/2}\dm{r})$ and, since they are positive,
they admit vector-valued extensions to $\ell^q$ too \citep[Theorem 5.5.10]{GrafakosCFA}.

There is one more operator to consider on p.\,280 of \cite{CordobaBalodis}:
\[
    T_{N,m}f(x) := \sqrt{x}f_{\nu}(A_N^\nu x)\int_0^1\sqrt{t}f_{\nu}(A_N^\nu t)f(t)\dm{t}
\]
with $f_\nu$ as in Lemma 1 of \cite{CordobaBalodis}.
For the stated range of $p$, $p/q\geq 1$, so we can apply Jensen's Inequality
as on p.\,280 of\cite{CordobaBalodis}, and we obtain the
$L^p_{\mathrm{rad}}(\ell_{\mathrm{ang}}^q)$-boundeness of $T_{N,\nu}$ too.

\begin{theorem}
    The operators $S_{N,M}$ are uniformly bounded on $L_{\mathrm{rad}}^p(\ell_{\mathrm{ang}}^q)$
    if $2\leq p<4$, when $N\geq AM+1$ for an absolute constant $A>0$.
\end{theorem}
\begin{cor}\label{cor:Lpconv}
    Let $2\leq p <4$. Then, if $N_k,M_k$ are sequences of natural numbers such that $N_k \geq AM_k+1$
    and $M_k\to\infty$, we have \[
        \lim_{k\to\infty}\norm{S_{N_k,M_k}f-f}_{p,q} = 0
    \] for all $f\in L_{\mathrm{rad}}^p(\ell_{\mathrm{ang}}^q)$.
    In particular, \[
        \lim_{k\to\infty}\norml{S_{N_k,M_k}f-f}{p}{\D} = 0
    \] for all $f\in L_{\mathrm{rad}}^p(\ell_{\mathrm{ang}}^q)$.
\end{cor}

Note, however, that this method does not allow us to conclude $L^p$ norm convergence \emph{for all} $f\in L^p$,
but only for the smaller space $L_{\mathrm{rad}}^p(\ell_{\mathrm{ang}}^q)\subsetneq L^p(\D)$.
The following function lies in $L^p$ \citep[since it is continuous, see][]{Zygmund2003}
but its Fourier coefficients are not $\ell^q$ summable for any $q>2$:
\[
    g(t) := \sum_{k=2}^\infty \frac{\e^{ik\log{k}}}{\sqrt{k}(\log{k})^2}\e^{\tpi t}.
\]In other words, $g\in L^p(\D)\setminus L_{\mathrm{rad}}^p(\ell_{\mathrm{ang}}^q)$.

The above example, which is a counterexample to Plancherel's theorem in $L^p$, hints at the underlying issue:
we did not obtain $L^p$ \emph{bounds} for the partial sum operators,
so we cannot apply the proof of \cref{cor:Lpconv} to $L^p(\D)$ directly.
\section{Concluding remarks}
To summarise, we have the following pieces of the convergence puzzle:

\begin{center}\begin{table}[ht]
    \begin{tabular}{c|cccccc}
        $p$ range & $[1,4/3)$ & $[4/3,2)$ & 2 & $(2,4)$ & 4&$(4,\infty)$\\
        \hline
        $\norml{\cdot}{p}{\D}$-convergence
        & No & ? & Yes & $f\in L_{\mathrm{rad}}^p(\ell_{\mathrm{ang}}^q)$ & ? & No
    \end{tabular}
    \caption{$L^p$ convergence of Bessel--Fourier series for various ranges of $p$.}
\end{table}\end{center}

In light of the existing results \citep{CordobaBalodis,Wing1950,BenedeckPanzone},
we might offer the following conjecture.
\begin{conj}
    For all $4/3<p<4$ and $f\in L^p(\D)$, \[
        \lim_{k\to\infty}\norml{S_{N_k,M_k}f-f}{p}{\D} = 0,
    \] for some appropriate choice of $N_k,M_k\in\N$, $k\geq 1$.
    For $p$ outside of this range, convergence fails in general.
\end{conj}
However, it is not completely clear that we can expect this.
\cite{Cordoba1989} proved that the ball multiplier is bounded on the mixed norm space
$L^p_{\mathrm{rad}}(L^2_{\mathrm{ang}})$ in the range $2n/(n+1)<p<2n/(n-1)$ that was originally conjectured for $L^p(\Rn)$,
but disproved in \cite{Fefferman1971}.
As we saw above, the $L^p_{\mathrm{rad}}(L^2_{\mathrm{ang}})$ argument works by essentially ``eliminating'' the angular eigenfunctions
and reducing the problem to bounds on the one-dimensional Bessel--Fourier series.
This ``trick'' is not available in $L^p(\D)$, owing to the failure of Plancherel's Theorem for $p\neq 2$,
so obtaining uniform bounds on the operators $S_{N,M}$ is considerably more difficult.

\section*{Acknowledgments}
I owe a special thanks to my doctoral advisor, Prof James C. Robinson,
for many insightful discussions and his helpful criticism of early drafts.
This work was supported by the EPSRC/2443915 studentship and the Warwick Mathematics Institute.

\bibliographystyle{unsrtnat}
\bibliography{bibliography}

\end{document}